\documentclass{amsart}

\date{}

\makeatletter

\RequirePackage{amssymb}

\newcommand{\theoremName}{Theorem}
\newcommand{\lemmaName}{Lemma}
\newcommand{\corollaryName}{Corollary}
\newcommand{\statementName}{Statement}
\newcommand{\remarkName}{Remark}
\newcommand{\exampleName}{Example}
\newcommand{\definitionName}{Definition}
\newcommand{\problemName}{Problem}
\newcommand{\proofName}{Proof}
\renewcommand{\proofname}{\proofName}
\newcommand{\answerName}{Answer}
\newcommand{\hintName}{Hint}

\theoremstyle{plain}

\newtheorem {theorem}{\theoremName}
\newtheorem {lemma}{\lemmaName}
\newtheorem {corollary}{\corollaryName}
\newtheorem {statement}{\statementName}

\theoremstyle{remark}

\newtheorem{Remark}{\remarkName}

\theoremstyle{definition}

\let\@newpf\proof \let\proof\relax 
\def \namepf[#1] {\@newpf[\proofname\ #1]}
\newenvironment{proof}{\@ifnextchar[{\namepf}{\@newpf[\proofname]}}{\qed\endtrivlist}

\def \Integer {{\mathbb Z}}

\def \Complex {{\mathbb C}}

\def \lnorm#1\rnorm {\vphantom{#1}\left\|\smash{#1}\right\|}
\def \lmod#1\rmod {\vphantom{#1}\left|\smash{#1}\right|}

\newcommand \bydef {\stackrel{\mbox{\scriptsize def}}{=}}

\@ifundefined{name}{\newcommand \name[1] {\mathop{\rm #1}\nolimits}}%
{\relax}

\newcommand \eps {\varepsilon}
\renewcommand \phi {\varphi}
\renewcommand \rho {\varrho}
\renewcommand \emptyset {\varnothing}

\def \Partit {{\mathcal D}}
\def \A {\mathcal A}
\def \plusmark {\text{``+''}}
\def \minusmark {\text{``--''}}
\newcommand{\set}[1]{\{1,2,\dots,#1\}}
\newcommand{\edge}[2]{\{#1,#2\}}
\numberwithin{equation}{section}

\makeatother

\begin{document}

\title{Around matrix-tree theorem}

\author[Yu.\,Burman]{Yurii Burman}
\address{Independent University of Moscow, Bolshoi  Vlas'evskii per. 11, 
121002, Moscow, Russia}
\email{burman@mccme.ru} 
\author[B.\,Shapiro]{Boris Shapiro}
\address{Mathematics Department, Stockholm University, S-106 91 
Stockholm, Sweden} 
\email{shapiro@math.su.se}

\thanks{Research supported in part by the RFBR grants \# N.Sh.1972.2003.1 
and \# 05-01-01012a}

 \begin{abstract}
Generalizing  the classical matrix-tree theorem we provide a formula 
counting subgraphs of a given graph with a fixed $2$-core. We use this 
generalization to obtain an analog of the matrix-tree theorem for the root 
system $D_n$ (the classical theorem corresponds to the $A_n$-case). Several 
byproducts of the developed technique, such as a new formula for a 
specialization of the multivariate Tutte polynomial, are of independent 
interest. 
 \end{abstract}

\maketitle

\section{Introduction}\label{Sec:Intro}

Let us first fix  some definitions and notation to be used throughout the 
paper. The main object of our study will be an undirected graph $G$ without 
multiple edges. It is understood as a subset $G \subset \{\{i,j\} \mid i,j 
\in \set{n}\}$, where elements of $\set{n}$ are vertices and elements of 
$G$ itself are edges. Informally speaking, this means that we mark (i.e.\ 
distinguish) vertices but not edges of $G$ (except for Section 
\ref{Sec:APoly} where an edge labeling  will be used). Usually we will 
assume that $G$ contains no loops, i.e.\ edges $\edge{i}{i}$. Directed 
graphs (appearing in Sections \ref{Sec:MSub} and \ref{Sec:Orient} for 
technical purposes) are subsets of $\set{n}^2$. Since a graph is understood 
as a set of edges, notation $F \subset G$ means that $F$ is a subgraph of 
$G$.

We will denote by $n = v(G)$ the number of vertices of $G$, by $\#G = e(G)$ 
the number of its edges, and by $k(G)$ the number of connected components. 
For every connected component $G_i \subset G$ ($i = 1, \dots, k(G)$) it 
will be useful to consider its Euler characteristics $\chi(G_i) = v(G_i) - 
e(G_i)$. A connected graph containing no cycles will be called a {\em 
tree}, a disconnected one, a {\em forest}. Note that the absence of cycles 
is equivalent to the equality $\chi(G_i) = 1$ for all $i$; if cycles are 
present then $\chi(G_i) \le 0$.

We will usually supply edges of the graph $G$ with weights. A weight 
$w_{ij} = w_{ji}$ of the edge $\edge{i}{j}$ is an element of any algebra 
$\A$. For a subgraph $F \subset G$ denote $w(F) \bydef \prod_{\edge{i}{j} 
\in F} w_{ij}$; call it the {\em weight} of $F$. For any set $U$ of 
subgraphs of $G$ call the expression $Z(U) = \sum_{F \in U} w(F)$ the {\em 
statistical sum} of $U$. (By definition, we assume $w_{ij} = 0$ if $G$ 
contains no edge $\edge{i}{j}$.)

To a graph $G$ with weighted edges one associates its {\em 
Laplacian matrix} $L_G$. It is a symmetric $(n \times n)$-matrix with the 
elements 
 \begin{equation*}
(L_G)_{ij} =  \begin{cases}
-w_{ij}, &i \ne j,\\
\sum_{k \ne i} w_{ik}, &i = j.
 \end{cases}
 \end{equation*}
The Laplacian matrix is degenerate; its kernel always contains the 
vector $(1,1,\dots,1)$. However, its principal minors are generally nonzero 
and enter  the classical matrix-tree theorem whose first version was proved 
by G.\,Kirchhoff in 1847:

 \begin{theorem}[\cite{Ki}] \label{Th:MTreeClass}
Let $T_G$ be the set of all (spanning) trees of $G$. Then $Z(T_G)$ is 
equal to any principal minor of $L_G$.
 \end{theorem}

This theorem has numerous generalizations (for a review, see e.g.\ 
\cite{Abde} and the references therein). For our purposes the most 
important will be the ``all-minors'' theorem by S.\,Chaiken \cite{Chaiken}. 

Call a subset $J = \{(i_1,j_1), \dots, (i_m,j_m)\} \subset 
\{1,2,\dots,n\}^2$ {\em component-disjoint} if $i_p \ne i_q$ and $j_p \ne 
j_q$ for every $p \ne q$; denote $\Sigma J \bydef \sum_{p=1}^m (i_p + 
j_p)$. Fix a numeration of the pairs $(i_p,j_p) \in J$ such that $i_1 < 
\dots < i_m$, and denote by $\tau_J$ the permutation of $\set{m}$ defined 
by the condition $j_{\tau_J(1)} < j_{\tau_J(2)} < \dots < j_{\tau_J(m)}$.

A forest $F$ with the vertex set $\set{n}$ is called {\em $J$-admissible} 
if it has $m$ components, and every component contains exactly one vertex 
from the set $\{i_1, \dots, i_m\}$, and exactly one, from $\{j_1, \dots, 
j_m\}$ (these two may coincide if the sets intersect). Denote by 
$\gamma_{F,J}$ a permutation of the set $\set{m}$ such that $i_p$ and 
$j_{\gamma_{F,J}(p)}$ lie in the same component of $F$, for every $p = 
1,2,\dots,m$. 

For an $n \times n$-matrix $M$ and a component-disjoint set $J$ denote by 
$M(J)$ the submatrix of $M$ obtained by deletion of the rows $i_1, \dots, 
i_m$ and the columns $j_1, \dots, j_m$. For any permutation $\sigma$ denote 
by $\eps(\sigma) = \pm1$ its sign (parity).

 \begin{theorem}[\cite{Chaiken}] \label{Th:AllMinors}
For any component-disjoint subset $J \subset \set{n}^2$ one has 
 \begin{equation*}
(-1)^{\Sigma J}\det (L_G)(J) = \sum_F \eps(\tau_J \circ \gamma_{F,J}) w(F)
 \end{equation*}
where the sum is taken over the set of all $J$-admissible subforests $F$ of 
$G$.
 \end{theorem}

Theorem \ref{Th:MTreeClass} is a particular case of Theorem 
\ref{Th:AllMinors} corresponding to the situation when $J$ contains one 
element only. 

Most of this article is devoted to various generalizations of Theorem 
\ref{Th:MTreeClass}. In Section \ref{Sec:MSub} we consider determinant-like 
expressions for statistical sums of subgraphs $F \subset G$ with cycles 
(namely, subgraphs with a given $2$-core). In Section \ref{Sec:ChiZero} we 
consider the case of subgraphs with vanishing Euler characteristics. 
Spanning trees of a graph $G$ can be interpreted as irreducible linearly 
independent subsets of roots in the root system $A_n$; in Section 
\ref{Sec:D} we prove an analog of Theorem \ref{Th:MTreeClass} for the root 
system $D_n$.

Two remaining sections form a sort of appendix to the paper. In Section 
\ref{Sec:Orient} we give an explicit formula for the number $d(G)$ of 
orientations of the graph $G$ without sources and sinks (this number enters 
Theorem \ref{Th:2Core}). In Section \ref{Sec:APoly} we provide a formula 
for the so called {\em external activity polynomial} which is a 
specialization of the multivariate Tutte polynomial of the graph $G$. The 
latter is defined as
 \begin{equation*}
T_G(q,w) = \sum_{m=1}^n q^n Z(U_m)
 \end{equation*}
(see \cite{Sokal,Tu,WelshMerino} for details) where $U_m$ is the set of 
subgraphs of $G$ having $m$ connected components and $w$ is a collection of 
weights of the edges. The formula we prove (Theorem \ref{Th:ExtAct}) is an 
alternating sign summation over partitions of the set of vertices of $G$. 

In the end of the paper we discuss several open problems related to the 
main topic.

\subsection*{Acknowledgments}

The first named author is sincerely grateful to the Mathematics 
Department of Stockholm University for the hospitality and financial 
support of his visit in September 2005 when the essential part of this 
project was carried out. We are thankful to Professors N.\,Alon and  
A.\,Sokal for their comments on the Tutte polynomial and a number of 
relevant references. We are grateful to Professor Olivier Bernardi who 
pointed out an important mistake in an earlier version of this paper.

\section{Graphs with a given $2$-core}\label{Sec:MSub}

Let $G$ be an undirected graph (loops and multiple edges are allowed). The 
maximal subgraph $G' \subset G$ such that every vertex of $G'$ is an 
endpoint of at least two edges or is attached to a loop (that is, there are 
no ``hanging'' vertices) is called the {\em $2$-core of $G$} and denoted by 
$\name{core}_2(G)$. A graph $G$ is the union of $\name{core}_2(G)$ and a 
number of forests (possibly empty) attached to every vertex of 
$\name{core}_2(G)$.

A graph $G$ is called {\em negative} if it contains no loops, no multiple 
edges, $G = \name{core}_2(G)$, and $\chi(G_i) < 0$ where $G_i$, $i = 1, 
\dots, k(G)$ are connected components of $G$. A graph $G$ is called {\em 
non-positive} if all the above is true but $\chi(G_i) \le 0$. 
A non-positive graph $G$ is the union of a negative graph $G_0$ and several 
cycles, each cycle forming a separate connected component. We will code this 
situation as  $G = G_0 \cup 3^{k_3} \dots n^{k_n}$ where $k_s$ stands for 
the total number of cycles of length $s$. 

For any directed graph $Q$ (with the vertex set $\set{n}$) denote by $[Q]$ 
the corresponding undirected graph. Given a $(n \times n)$-matrix $M$ with 
entries $a_{i,j}\in \A$ define
 \begin{equation*}
\langle M,Q\rangle \bydef \prod_{(i,j) \in Q} a_{ij}.
 \end{equation*}
In particular, if $[Q] \subset G$ where $G$ is a graph without loops or 
multiple edges, with weights $w_{ij}$ (like in the previous section), then 
$\langle L_G,Q\rangle = (-1)^{e(Q)}w(Q)$.

A directed graph $Q$ is called {\em regular} if the following two 
conditions are satisfied:

 \begin{enumerate}
\item\label{It:NSS} $Q$ contains no sources or sinks, i.e.\ for every 
vertex there is at least one incoming  and one  outgoing edge.

\item If $Q$ contains a loop (an edge $(i,i)$) or a pair of antiparallel 
edges (edges $(i,j)$ and $(j,i)$) then they form a separate connected 
component of $Q$.
 \end{enumerate}
If $Q$ is a regular directed graph then $[Q]$ consists of a non-positive 
graph and several loops and double edges (cycles of length $2$), each loop 
and double edge forming a separate connected component. We will denote this 
by $[Q] = H \cup 1^{k_1} 2^{k_2}$ where $H$ is non-positive and $k_1, k_2$ 
are the number of loops and double edges, respectively.

In what follows it will be convenient to allow graphs to have multiple 
(more specifically, double) edges. If $F$ is a graph with multiple edges we 
will abuse notation writing $F \subset G$ if the graph obtained from $F$ by 
neglecting the multiplicities is a subgraph of $G$. Computing the weights, 
we will, however, take multiplicities into account: 
 \begin{equation*}
w(F) \bydef \prod_{\edge{i}{j} \text{ is an edge of }F} w_{ij}^{m_{ij}}
 \end{equation*}
where $m_{ij} \in \Integer_{{\ge}0}$ is the multiplicity of the edge 
$\edge{i}{j}$.

Let $H \subset G$ be a non-positive graph plus several double 
edges, each double edge forming a separate component. In other words, 
$H = H_0 \cup 2^{k_2} 3^{k_3} \dots n^{k_n}$ where $H_0$ is negative. Then 
denote
 \begin{equation*}
\rho_G(H) = \sum_{\scriptsize\begin{array}{c}
\Lambda \mbox{ is regular}\\
\relax[\Lambda] = H \cup 1^{n-v(H)}
\end{array}}
\langle L_G,\Lambda\rangle
 \end{equation*}
(so that the total number of vertices of $\Lambda$ is $n$).  By $U(H)$ 
denote the set of all subgraphs $F \subset G$ such that $\name{core}_2(F) = 
H$. 

 \begin{theorem} \label{Th:2Core}
Let $H = H_0 \cup 2^{k_2} 3^{k_3} \dots n^{k_n}$ be a non-positive graph 
without loops together with several double edges. Then 
 \begin{equation}\label{Eq:MSub}
\rho_G(H) = (-1)^{e(H)} d(H_0) \sum_{l_2=k_2}^n \dots \sum_{l_n=k_n}^n 
\binom{l_2}{k_2} \dots \binom{l_n}{k_n} 2^{l_3 + \dots + l_n} \, 
Z(U(H_0 \cup 2^{l_2} \dots n^{l_n}))
 \end{equation}
where $d(H_0)$ is the number of orientations of $H_0$ without sources and 
sinks. 
 \end{theorem}

 \begin{corollary} \label{Cr:MSubgrInv}
One has 
 \begin{equation}\label{Eq:MSubgrInv}
 \begin{aligned}
Z(U(H))  =(-1)^{e(H_0)} d(H_0) 2^{-(k_3 + \dots + k_n)} &\times\\
\times \sum_{l_2=k_2}^n \dots \sum_{l_n=k_n}^n (-1)^{l_2 + l_3 + \dots + l_n}& 
\binom{l_2}{k_2} \dots \binom{l_n}{k_n} \, \rho_G(H_0 \cup 2^{l_2} 
\dots n^{l_n}).
 \end{aligned}
 \end{equation}
 \end{corollary}

 \begin{Remark}
Corollary \ref{Cr:MSubgrInv} is our closest approximation  to a 
``matrix-subgraph'' theorem, that is, the best available analog of Theorem 
\ref{Th:MTreeClass} for subgraphs of arbitrary structure. Indeed, the 
left-hand side of \eqref{Eq:MSubgrInv} is the statistical sum over the 
graphs with a fixed $2$-core (for trees the $2$-core is empty), while the 
right-hand side is a polylinear function of matrix elements of the 
Laplacian matrix (in the case of trees it was its principal minor). Notice 
that, unlike Theorem \ref{Th:MTreeClass}, the right-hand side of 
\eqref{Eq:MSubgrInv} cannot be computed in polynomial time. This is hardly 
surprising: it is known that the calculation of the Tutte polynomial (and 
even its value at almost any point of the plane) is a sharp $P$-hard 
problem (see \cite[\S9]{Welsh}). Therefore there is no hope to obtain a 
formula for the statistical sum of connected subgraphs in $G$ with any 
given number of edges in the form of a determinant or, in general, to get a 
formula of polynomial complexity in terms of the Laplacian matrix. 
 \end{Remark}

 \begin{proof}[of Theorem \ref{Th:2Core}]
Let $\Lambda = \Lambda_0 \cup \Lambda_1$ be a regular subgraph of $G$ such 
that $[\Lambda_0] = H$ and $[\Lambda_1] = 1^{n-v(H)}$. Now, $\langle 
L_G,\Lambda\rangle = \langle L_G,\Lambda_0\rangle \langle 
L_G,\Lambda_1\rangle$. Since $\Lambda_0$ contains no loops, then $\langle 
L_G,\Lambda_0\rangle = (-1)^{e(H)} w(H)$.

One has $(L_G)_{ii} = \sum_{k \ne i} w_{ik}$, so that the term $\langle 
L_G,\Lambda_1\rangle$ can be represented as the sum of monomials $w_{i_1k_1} 
\dots w_{i_s k_s}$ where $\{i_1, \dots, i_s\}$ is the vertex set of 
$\Lambda_1$. In other words, $\langle L_G,\Lambda_1\rangle = \sum_{\Theta} 
w(\Theta)$ where $\Theta$ is the  directed graph with $[\Theta] \subset G$ 
satisfying the following property: if $i \in \{i_1, \dots, i_s\}$ then 
$\Theta$ contains exactly one edge starting from $i$, and if $i \notin 
\{i_1, \dots, i_s\}$ is a vertex of $\Theta$ then it is a sink (no edge 
starts from it). 

One can easily see that every connected component of $\Theta$ is either a 
tree such that all its vertices except the root are in $\{i_1, \dots, 
i_s\}$, or a graph with exactly one cycle with all its vertices in $\{i_1, 
\dots, i_s\}$. Thus, $\name{core}_2([\Lambda_0 \cup \Theta]) = H_0 \cup 
2^{l_2} \dots n^{l_n},$ where $l_2 \ge k_2$, \dots, $l_n \ge k_n$.

On the other hand, let $F \subset G$ be a subgraph such that 
$\name{core}_2(F) = H_0 \cup 2^{l_2} \dots n^{l_n}$. To identify $F$ with 
$[\Lambda_0 \cup \Theta]$ one has, first, to point out which ``$1$-cycled'' 
connected components of $F$ belong to $\Lambda_0$ and which to $\Theta$ --- 
there are $\binom{l_2}{k_2} \dots \binom{l_n}{k_n}$ ways to do this. Having 
this choice made one must orient the $2$-core of $F$ without sources and 
sinks --- the number of such orientations being $d(H_0 \cup 2^{l_2} \dots 
n^{l_n}) = 2^{l_3 + \dots + l_n} d(H_0)$.
 \end{proof}

 \begin{proof}[of Corollary \ref{Cr:MSubgrInv}]
One has $\frac{1}{k!}\frac{d^k}{d x^k} (1-x)^l = \sum_{s=0}^l \binom{s}{k} 
\binom{l}{s} (-1)^s x^s = \binom{l}{k} (1-x)^{l-k}$, and therefore 
 \begin{equation*}
\sum_{s=0}^l (-1)^s\binom{s}{k} \binom{l}{s} =  \begin{cases}
0, &\mbox{if } k < l,\\
1, &\mbox{if } k = l.
 \end{cases}
 \end{equation*}
The corollary is now straightforward.
 \end{proof}

\section{Graphs with vanishing Euler characteristics}\label{Sec:ChiZero}

Corollary \ref{Cr:MSubgrInv} becomes particularly simple if $H$ is a cycle. 
Namely, if $H = s^1$ (a cycle of length $s$) then $\lambda_G(H)$ is the 
statistical sum of the set of all connected subgraphs $F \subset G$ having 
exactly one cycle of length $s$. The ``negative part'' $H_0$ of the graph $H$ 
is empty which implies $d(H_0) = 1$. 

Denote by $\Sigma_n$ the symmetric group of order $n$ acting on $\set{n}$, 
and denote by $\Partit_n$ the set of all partitions of $n$. For a 
permutation $\sigma \in \Sigma_n$ having $k_1$ cycles of length $1$, $k_2$ 
cycles of length $2$, etc., denote $D(\sigma) \bydef 1^{k_1} \dots n^{k_n} 
\in \Partit_n$. Finally, for any function $f: \Partit \to \A$ define the 
{\em $f$-determinant} of an $(n \times n)$-matrix $M$ with entries 
$a_{i,j}\in \A$ by the formula 
 \begin{equation*}
\det\nolimits_f (L) = \sum_{\sigma \in \Sigma_n} f(D(\sigma)) 
a_{1,\sigma(1)} \dots a_{n,\sigma(n)}.
 \end{equation*}
Now one has
 \begin{align*}
\rho_G(2^{l_2} \dots n^{l_n}) &= \sum_{D(\sigma) = 1^{n - 2l_2 - 
\dots - nl_n}2^{l_2} \dots n^{l_n}} (-1)^{n + 2l_2 + \dots + nl_n} 
(L_G)_{1,\sigma(1)} \dots (L_G)_{n,\sigma(n)} \\
&= \det\nolimits_{\chi_{l_2, \dots, l_n}} L_G,
 \end{align*}
where 
 \begin{equation*}
\chi_{l_2, \dots, l_n}(1^{k_1} \dots n^{k_n}) = 
 \begin{cases} 
(-1)^{n + 2l_2 + \dots + nl_n},  &\text{if } k_2 = l_2, k_3=l_3, \dots, k_n 
= l_n,\\ 
0, &\text{otherwise.}
 \end{cases}
 \end{equation*}
Thus, Corollary \ref{Cr:MSubgrInv} for a cycle takes the following form: 

 \begin{statement} \label{St:1Cycle}
The statistical sum of the set of subgraphs $F \subset G$ having one 
cycle of length $s \ge 3$ is equal to $\frac12 \det\nolimits_{\tau_s} L_G$, 
where $\tau_s(1^{k_1} \dots n^{k_n}) = (-1)^{n + 2k_2 + \dots + nk_n} k_s$. 
The statistical sum of the set of subgraphs $F \subset G$ having one 
cycle of length $2$ is $\det\nolimits_{\tau_2} L_G$.
 \end{statement}

This corollary implies the following formula which is the 
``matrix-tree theorem'' for connected subgraphs containing exactly one cycle 
of any length $s \ge 3$, that is, connected subgraphs $H \subset G$ with 
$\chi(H) = 0$: 

 \begin{corollary} 
Let $U_G$ be the set of all connected subgraphs of $H \subset G$ such that 
$\chi(H) = 0$. Then 
 \begin{equation*}
Z(U_G) = \frac{1}{2}\det\nolimits_\mu(L_G) 
  \end{equation*}
where 
 \begin{equation}\label{Eq:1Cycle}
\mu(1^{k_1} 2^{k_2} \dots n^{k_n}) = (-1)^{n + k_2 + 2k_3 + \dots + (n-1)k_n} 
(2k_2 + k_3 + \dots + k_n).
 \end{equation}
 \end{corollary}

A finer result concerning graphs $H \subset G$ such that $\chi(H_i) = 0$ 
for any connected component $H_i$ of $H$ ($i = 1, \dots, k(H)$) can be 
obtained using Theorem \ref{Th:AllMinors}. 

For a graph $G$ and a component-disjoint set $J$ denote by $G - J$ the 
graph obtained from $G$ by deletion of all the edges $(i_p j_p)$ where 
$(i_p,j_p) \in J$. Then Theorem \ref{Th:AllMinors} implies

 \begin{statement} \label{St:GivenCycle}
Let $G$ be a graph with the vertex set $\set{n}$, without loops and 
multiple edges, with weights $w_{ij}$ defined for all the edges. Let $J = 
\{(i_1,j_1), \dots, (i_m,j_m)\}$ be a component-disjoint subset of 
$\set{n}^2$. Then
 \begin{equation}\label{Eq:GivenCycle}
(-1)^n \eps(\tau_J) w_{i_1 j_1} \dots w_{i_m j_m} \det (L_{G-J})(J) = 
\sum_H (-1)^{k(H)}w(H)
 \end{equation}
where the sum is taken over the set of all subgraphs $H\subseteq G$ such 
that every connected component $H_i$ of $H$ contains one cycle
(that is, $\chi(H_i) = 0$), the edges $\edge{i_1}{j_1}, \dots, 
\edge{i_m}{j_m}$ enter these cycles and vertices $i_p$ and $j_q$ alternate 
along the cycle.
 \end{statement}

 \begin{proof}
It follows from Theorem \ref{Th:AllMinors} that the product $w_{i_1 j_1} 
\dots w_{i_m j_m} \det (L_{G-J})(J)$ is equal to the sum of $\pm w_{i_1 
j_1} \dots w_{i_m j_m} w(F)$ where $F$ runs over the set of subforests of 
$G - J$ having $m$ components and such that the $p$-th component contains 
the vertices $i_p$ and  $j_{\gamma_{F,J}(p)}$; here $\gamma_{F,J}$ is the 
permutation of $\{1,2,\dots,m\}$ defined in Section \ref{Sec:Intro}. In 
other words, $w_{i_1 j_1} \dots w_{i_m j_m} \det (L_{G-J})(J)$ is equal to 
the sum of $\pm w(H)$ where $H = F+J$ is the result of addition to $F$ of 
the edges $\edge{i_1}{j_1}$, \dots, $\edge{i_m}{j_m}$. Thus, $H$ is a graph 
with one cycle in every connected component; all edges $\edge{i_p}{j_p}$ 
enter the cycles, and vertices $i_p$ and $j_q$ alternate 
along the cycle. The connected components of $H$ are in one-to-one 
correspondence with the cycles of the permutation $\gamma_{F,J}$. The sign 
of the term $w(F)$ is equal to $(-1)^n \eps(\tau_J) \eps(\tau_J \circ 
\gamma_{F,J}) = (-1)^n \eps(\gamma_{G,J})$. The permutation $\gamma_{G,J}$ 
contains $k(H)$ cycles. The sign of any permutation of $\set{n}$ with $k$ 
cycles equals $(-1)^{n+k}$, and therefore, the total sign is $(-1)^{k(H)}$.
 \end{proof}

Denote now
 \begin{equation*}
Q_m = \sum_{\#J = m} w_{i_1 j_1} \dots w_{i_m j_m} \det (L_{G-J})(J)
 \end{equation*}
where the sum is taken over the set of all component-disjoint subsets $J 
\subset \set{n}^2$ of cardinality $m$. Statement \ref{St:GivenCycle} allows 
to express the generating function for the sequence $Q_m$:

 \begin{theorem} \label{Th:GenFun}
One has 
 \begin{equation}\label{Eq:GenFun}
\sum_{m=1}^\infty Q_m t^m = (-1)^n \sum_H w(H) \prod_{i=1}^{k(H)} 
((1+t)^{\ell_i(H)}-1)
 \end{equation}
where the sum in the right-hand side is taken over the set of all subgraphs 
$H \subset G$ such that $\name{core}_2(H_i)$ is a cycle of 
length $l_i(H)$; here $H_1, \dots, H_{k(H)}$ are connected components of 
$H$.
 \end{theorem}

 \begin{proof}
By Statement \ref{St:GivenCycle} one has that $Q_m = \sum_H a_m(H) w(H)$ 
where the sum is taken over the set of all subgraphs $H \subset G$ having 
exactly one cycle in every connected component. The coefficient $a_m(H)$ is 
equal, to $(-1)^{n+k(H)}$ times the number of component-disjoint sets $J = 
\{(i_1,j_1), \dots, (i_m,j_m)\}$ such that 
 \begin{itemize}
\item $\edge{i_p}{j_p} \in \name{core}_2(H)$ for all $p = 1, \dots, m$.

\item For every cycle of $H$ there is at least one edge $(i_p j_p)$ 
entering it.

\item If a cycle of $H$ has more than one edge $(i_p, j_p)$ in it then the 
vertices $i_p$ and  $j_q$ alternate  along the cycle.
 \end{itemize}
This obviously implies that 
 \begin{equation*}
a_m(H) = (-1)^{n+k(H)} 
\sum_{\scriptsize\begin{array}{l} 
m_1 + \dots + m_{k(H)} = m\\ 
m_1, \dots, m_{k(H)} \ge 1
 \end{array}}\binom{\ell_1(H)}{m_1} \dots 
\binom{\ell_{k(H)}(H)}{m_{k(H)}},
 \end{equation*}
and \eqref{Eq:GenFun} follows.
 \end{proof}

 \begin{corollary}
 \begin{equation}\label{Eq:AltSum}
\sum_{m=1}^\infty (-1)^m Q_m = \sum_H (-1)^{n+k(H)} w(H).
 \end{equation}
 \end{corollary}

\section{Linearly independent subsets of the root systems $A_n$ and 
$D_n$}\label{Sec:D}

The technique of Section \ref{Sec:ChiZero} can be used to obtain  results 
on linearly independent subsets of finite root systems, cf.\ 
\cite{PoShSh}.  

The set of positive roots $R_+(A_n)$ of the reflection group $A_n$ consists 
of vectors $e_{ij} = b_i - b_j$, $1 \le i < j \le n$ where $b_1, \dots, 
b_n$ is the standard basis in $\Complex^n$. We will assign to every root 
$e_{ij} \in R_+(A_n)$ its weight $w_{ij} \in \A$ where $\A$ is any algebra. 
By definition $w_{ji} = w_{ij}$. For any subset $S \subset R_+(A_n)$ of 
positive roots consider a graph $\Gamma(S)$ with the vertices $1, 
\dots, n$ such that $\edge{i}{j}$ is an edge of $\Gamma(S)$ wherever 
$e_{ij} \in S$. The edge $\edge{i}{j}$ bears the weight $w_{ij}$. The graph 
$\Gamma(S)$ is undirected and contains no loops or multiple edges. If $S' 
\subset S$ then $\Gamma(S')$ is a subgraph of $\Gamma(S)$. We will write 
$w(S)$ instead of $w(\Gamma(S))$ for short and denote by $L_S$ the 
Laplacian matrix of the graph $\Gamma(S)$. 

For a given subset $S\subset R_+(A_n)$ one can consider the group $G(S)$ 
generated by the reflections in the roots $e_{ij} \in S$. The group $G(S)$ 
is a subgroup of the Weyl group of $A_n$, and therefore the space $V = 
\{\sum_{i=1}^n x_i b_i \mid \sum_{i=1}^n x_i = 0\} \subset \Complex^n$ is 
$G(S)$-invariant. $S$ is called {\em irreducible} if $V$ is an irreducible 
representation of $G(S)$.

The following is obvious:

 \begin{theorem} \label{Th:A}
A set $S \subset R_+(A_n)$ is linearly independent if and only if 
$\Gamma(S)$ contains no cycles. $S$ is irreducible if and only if 
$\Gamma(S)$ is connected. A linearly independent set $S' \subset S$ is 
maximal (among linearly independent subsets of $S$) if and only if 
$\Gamma(S')$ is a forest composed of spanning trees of connected components 
of $\Gamma(S)$. If $S$ is irreducible (that is, $\Gamma(S)$ connected) then 
any maximal linearly independent subset $S'$ of $S$ is also irreducible 
(that is, $\Gamma(S')$ is a spanning tree of $\Gamma(S)$).
 \end{theorem}

Using matroid terminology, one can reformulate Theorem \ref{Th:A} as 
follows. (See \cite{Sokal,WelshMerino} for more detail about matroids.) 

 \begin{corollary}
A submatroid of the linear matroid of $\Complex^n$ generated by vectors 
$e_{ij} \in S$ is isomorphic to the graphical matroid of $\Gamma(S)$. 
 \end{corollary}

One can associate a weight $w_{ij} = w_{ji} \in \A$ to every
root $e_{ij} \in R_+(A_n)$. So, one can consider weights of the root 
systems and statistical sums of sets of root systems, like it was
done for graphs in the previous sections. Now the matrix-tree theorem (i.e.\ 
Theorem~\ref{Th:MTreeClass}) and Theorem~\ref{Th:A} imply: 

 \begin{statement}\label{St:NTrees}
Let $S \subset R_+(A_n)$ be irreducible and $T_S$ be the collection of all 
maximal linearly independent subsets of $S$. Then $Z(T_S)$ is equal to 
(any) principal minor of the Laplacian matrix $L_S$.
 \end{statement}

Consider now a similar question for the reflection group  $D_n$. Its set 
$R_+(D_n)$ of positive roots consists of the vectors $e_{ij}^+ = b_i - b_j$ 
(the $\plusmark$-vectors) and $e_{ij}^- = b_i + b_j$ (the 
$\minusmark$-vectors) for all $1 \le i < j \le n$. We associate to every 
$\plusmark$-vector $e_{ij}^+$ the weight $u_{ij} \in \A$, and to every 
$\minusmark$-vector $e_{ij}^-$ the weight $v_{ij} \in \A$. Notions of  
linearly independent, maximal and irreducible subsets $S \subset R_+(D_n)$ 
are defined exactly as in the $A_n$-case.

For every set $S \subset R_+(D_n)$ consider the graph $\Gamma(S)$ with the 
vertices $1, \dots, n$ where the vertices $i$ and $j$ are joined by the 
edge marked $\plusmark$ if $e_{ij}^+ \in S$, and by the edge marked 
$\minusmark$ if $e_{ij}^- \in S$. Thus, the graph $\Gamma(S)$ is 
undirected, contains no loops, and has at most two edges joining every pair 
of vertices; all its edges are marked by $\plusmark$ or $\minusmark$, and 
if two edges join the same pair of vertices then their marks are different. 

A cycle in $\Gamma(S)$ is called {\em odd} if it contains an odd number of 
edges marked $\minusmark$. 

 \begin{theorem} \label{Th:D}
A set $S \subset R_+(D_n)$ is irreducible if and only if $\Gamma(S)$ is 
connected. $S$ is linearly independent if and only if every connected 
component of $\Gamma(S)$ is either a tree or a graph with exactly 
one cycle, and this cycle is odd. If $S$ is irreducible then a linearly 
independent set $S' \subset S$ is maximal if and only if the following 
holds: if $\Gamma(S)$ contains no odd cycles then $S' = S$, otherwise 
every connected component of $\Gamma(S')$ is a graph containing exactly one 
cycle, and this cycle is odd.
 \end{theorem}

This is a $D$-analog of Theorem \ref{Th:A} and it is obvious as well. Our goal in 
this section is to obtain a $D$-analog of Statement \ref{St:NTrees}.

Let $J = \{(i_1, j_1), \dots, (i_m,j_m)\} \subset \set{n}^2$ be a 
component-disjoint subset. Denote $S-J^- \bydef S \setminus \{e^-_{i_1 
j_1}, \dots, e^-_{i_m j_m}\}$.

 \begin{theorem} \label{Th:Cardm}
One has 
 \begin{equation*}
\sum_{m=1}^n t^m \sum_{\#J = m} v_{i_1 j_1} \dots v_{i_m j_m} \det (L_{S - 
J^-})(J) = (-1)^n \sum_H w(H) \prod_{i=1}^{k(H)} 
((1+t)^{\ell_i^-(H)}-1). 
 \end{equation*}
Here the internal sum in the left-hand side is taken over the set of all 
component-disjoint sets $J \subset \set{n}^2$ of cardinality $m$.  The sum 
in the right-hand side is taken over the set of all subsets $H \subset S$ 
such that every connected component of the graph $\Gamma(H)$ contains 
exactly one cycle. Above we denote by $k(H) \bydef k(\Gamma(H))$ the number 
of these components, and by $\ell_i^-(H)$ ($i = 1, \dots, k(H)$) the number 
of $\minusmark$-edges entering the cycle in the $i$-th component. 
 \end{theorem}

The proof is completely analogous to that of Theorem \ref{Th:GenFun}. A 
required analog of Statement \ref{St:NTrees} is now:

 \begin{corollary}
 \begin{equation}\label{Eq:SumD}
\sum_{m=1}^n (-2)^m \sum_{\#J = m} v_{i_1 j_1} \dots v_{i_m j_m} \det (L_{S 
- J^-})(J) = (-1)^n \sum_F (-2)^{k(F)} w(F).
 \end{equation}
Here the internal sum in the left-hand side is taken over the set of all 
component-disjoint sets $J \subset \set{n}^2$ of cardinality $m$. The sum 
in the right-hand side is taken over the set of all maximal linearly 
independent subsets $F \subset S$. As usual, $k(F)$ is the number 
of connected components of the graph $\Gamma(F)$.
 \end{corollary}

 \begin{proof}
This follows directly from Theorems \ref{Th:Cardm} and \ref{Th:D} and the 
equality
 \begin{equation*}
\prod_{i=1}^k ((-1)^{\ell_i}-1) =  \begin{cases}
(-2)^k, &\text{ if all the $\ell_i$ are odd,}\\
0, &\text{ if at least one $\ell_i$ is even.}
 \end{cases}
 \end{equation*}
 \end{proof}

\section{Orientations without sources and sinks}\label{Sec:Orient}

Let $G$ be an undirected graph with the vertex set $\set{n}$, without loops 
(multiple edges are allowed). In this section we give a combinatorial 
description of the number $d(G)$ of directed graphs $Q$ such that $[Q] = G$ 
and $Q$ has no sources or sinks. (Recall that the number $d(G)$ enters 
equation \eqref{Eq:MSub}.)

For a set of vertices $P \subset \set{n}$ of $G$ denote by $\langle 
P\rangle$ the subgraph of $G$ spanned by $P$ (i.e.\ having $P$ as its 
vertex set and containing all the edges of $G$ with both endpoints in $P$). 
Denote $k(P) \bydef k(\langle P\rangle)$ for short and denote by 
$\mu(P) = e(\langle \set{n} \setminus P\rangle)$, that is, the number of 
edges in $G$ having both endpoints outside $P$. 

Recall that a graph $F$ is called {\em bipartite} if  one can split its 
vertices into two groups such that every edge joins two vertices from 
different groups. Equivalently, this means that every closed path in $F$ 
contains an even number of edges.

 \begin{theorem}
Assume that $G$ has no isolated vertices. Then the number $d(G)$ of 
orientations of $G$ without sources and sinks (i.e.\ such that for every 
vertex there is at least one incoming and one outgoing edge) is given by 
the expression
 \begin{equation*}
d(G) = \sum_{m=0}^n (-1)^m \sum_{\scriptsize P: \#P = m, \langle P\rangle 
\mbox{ is bipartite}} 2^{\mu(P) + k(P)}.
 \end{equation*}
(By assumption, $k(\emptyset) = 0$ and $\mu(\emptyset)$ is equal to the 
total number of edges in $G$.)
 \end{theorem}

 \begin{proof}
Fix a set $P$ of vertices, and let $N(P)$ be the number of orientations of 
$G$ such that every vertex from $P$ is either a source or a sink. Since an 
edge cannot join two sources or two sinks, one has $N(P) = 0$ if $\langle 
P\rangle$ is not bipartite.

Suppose now that $\langle P\rangle$ is bipartite. Consider the graph $F$ 
obtained by adding to $\langle P\rangle$ all the edges having one vertex in 
$P$ and the other outside $P$. Apparently, $k(F) = k(P)$. Since $G$ has no 
isolated vertices, every connected component of $F$ has $2$ orientations 
such that every its vertex is either a source or a sink. Thus, the total 
number  of ways to orient the edges of $F$ is $2^{k(P)}$. The number of 
edges of $G$ not belonging to $F$ is $\mu(P)$. These edges can be oriented 
arbitrarily, and, therefore, $N(P) = 2^{\mu(P) + k(P)}$. The statement 
follows now from the inclusion-exclusion formula.
 \end{proof}

 \begin{corollary} \label{Cr:ChrPoly}
The number of orientations of $G$ without sources and sinks is given by 
 \begin{equation}\label{Eq:TotCyclChr}
d(G) = \sum_{m=0}^n (-1)^m \sum_{P: \#P = m} 2^{\mu(P)} \name{chr}_{\langle 
P\rangle}(2)
 \end{equation}
where $\name{chr}_F(\lambda)$ is the chromatic polynomial of the graph $F$, 
that is, the number of ways to color its vertices in $\lambda$ colors so 
that any two adjacent vertices have different colors. (One assumes 
$\name{chr}_{\langle \emptyset\rangle} = 1$.) 
 \end{corollary}

 \begin{proof} 
One has $\name{chr}_F(2) = 2^{k(F)}$ if the graph $F$ is bipartite, and 
$\name{chr}_F(2) = 0$ otherwise.
 \end{proof}

 \begin{corollary}
The number of orientations of $G$ without sources and sinks is given by 
 \begin{equation}\label{Eq:TotCyclChi}
d(G) = \sum_{F \subseteq G} 2^{\mu(F) + k(F)} (-1)^{\chi(F)}
 \end{equation}
where the sum is taken over the set of all subgraphs $F \subseteq G$, and 
$\mu(F)$ is the total number of edges in $G$ having no common vertices  
with the edges from $F$.
 \end{corollary}

 \begin{proof}
A classical result (see e.g.\ \cite{WelshMerino} for proof) relates the 
multivariate Tutte polynomial to the chromatic polynomial:
 \begin{equation*}
\name{chr}_H(\lambda) = T_H(\lambda,-1) = \sum_{\scriptstyle 
 \begin{array}{c}
\scriptstyle F \subseteq H\\
\scriptstyle v(F) = v(H)
 \end{array}} \lambda^{k(F)} (-1)^{e(F)}.
 \end{equation*}
($-1$ in the argument of $T_H$ means that one takes $w_{ij} = -1$ for every 
edge $\edge{i}{j}$ of $H$). Now by Corollary \ref{Cr:ChrPoly}, 
 \begin{equation*}
d(G) = \sum_{m = 0}^n (-1)^m \sum_{\scriptstyle \begin{array}{c}
\scriptstyle F \subseteq G\\
\scriptstyle v(F) = m
\end{array}} 2^{\mu(F)} 2^{k(F)} (-1)^{e(F)} = \sum_{F \subseteq G} 
2^{\mu(F)+k(F)} (-1)^{\chi(F)}.
 \end{equation*}
where $\chi(F) = v(F) - e(F)$ is the Euler characteristics of $F$.
 \end{proof}

\section{Multivariable external activity polynomial}\label{Sec:APoly}

In the previous section we made use of the fact that the chromatic 
polynomial of a graph is a specialization of its multivariate Tutte 
polynomial $T(q,w)$. Below we consider another specialization of the $T$ 
which we call the {\em external activity polynomial}.

As in Section \ref{Sec:MSub} let $G$ be a graph without loops or multiple 
edges and with the weights $w_{ij} \in \A$ assigned to its edges. 
Suppose that $G$ is connected and fix an arbitrary 
numeration of its edges. Let $T$ be a spanning tree of $G$ and $e$ be an 
edge not entering $T$. The graph $T \cup e$ has exactly one cycle, 
and this cycle contains the edge $e$. An edge $e$ is called {\em externally 
active} for $T$ if $e$ is the smallest edge (with respect to the above 
numeration) in the cycle. The polynomial
 \begin{equation*}
C_G(w) = \sum_{T \text{ is a spanning tree of } G} w(T) 
\prod_{\edge{i}{j} \text{ is externally active for } T} (w_{ij}+1).
 \end{equation*}
will be called the {\em external activity polynomial} of $G$. Its 
specialization (re)appeared recently in the form of the Hilbert polynomial 
for a certain commutative algebra related to $G$, see \cite{PoSh}.

Obviously, the following statement holds:

 \begin{statement} \label{St:Connect}
If $G$ is connected then $C_G(w) = Z(U_0)$, where $U_0$ is the set of all
connected spanning subgraphs of $G$.
 \end{statement}

 \begin{corollary}[\cite{Sokal}] \label{Cr:FreeTerm}
If $G$ is connected then $C_G(w) = \lim_{q \to 0} T_G(q,w)/q$ where $T_G$ 
is the multivariate Tutte polynomial.
 \end{corollary}

We present  another expression for the polynomial $C_G$:

 \begin{theorem}\label{Th:ExtAct} 
 \begin{equation}\label{Eq:ExtAct} 
C_G(w) = \sum_{k=1}^n (-1)^{k-1}(k-1)! \sum_{\set{n} = P_1 \sqcup 
\dots \sqcup P_k}\, \prod_{\edge{i}{j}: \exists s\,\, i,j \in P_s} (w_{ij}+1)
 \end{equation}
where the internal sum taken is over all partitions of the set of 
vertices into $k \ge 1$ pairwise disjoint subsets, and the product is taken 
over the set of all edges $\edge{i}{j}$ of the graph $G$ such that both 
endpoints ($i$ and $j$) belong to some $P_s,\,s = 1, \dots, k$.
 \end{theorem}

To prove Theorem \ref{Th:ExtAct} we need the following technical lemma.

 \begin{lemma}\label{Lm:Moebius}
For any $n \ge 2$ one has
 \begin{equation}\label{Eq:Moebius}
\sum_{k=1}^n (-1)^{k-1} \sum_{\set{n} = P_1 \sqcup \dots \sqcup P_k} (p_1 - 
1)! \dots (p_k - 1)! = 0, 
 \end{equation}
where $p_i = \# P_i$ is the cardinality of $P_i$. 
 \end{lemma}

 \begin{proof}
Denote by $a_{k,n}$ the coefficient at $(-1)^{k-1} (k-1)!$ in 
\eqref{Eq:Moebius} and use induction on $n$ to prove the lemma. For $n = 2$ 
one has $a_{1,2} = 1$ (the only  possible set partition is  
$\{1,2\} = \{1,2\}$) and $a_{2,2} = 1$ (the only possible set 
partition is $\{1,2\} = \{1\} \sqcup \{2\}$), so that \eqref{Eq:Moebius} 
holds. Let now it hold for some $n$. Partitions of the set $\set{n+1}$ fall 
into two types: either $\set{n+1} = P_1 \sqcup \dots \sqcup P_k \sqcup 
\{n+1\}$, or $\set{n+1} = P_1 \sqcup \dots \sqcup (P_s \cup \{n+1\}) \sqcup 
\dots \sqcup P_k$, where in both cases $P_1 \sqcup \dots \sqcup P_k$ is a 
partition of $\set{n}$. The sum taken over the set partitions of 
the first type is 
 \begin{equation*}
\sum_{k=1}^n (-1)^k \sum_{\set{n+1} = P_1 \sqcup \dots \sqcup P_k \sqcup 
\{n+1\}} (p_1 - 1)! \dots (p_k - 1)!\, 0! = 0
 \end{equation*}
by the induction hypothesis. The sum taken over the set partitions 
of the second type equals 
 \begin{align*}
\sum_{k=1}^n &(-1)^{k-1} \sum_{s=1}^k \sum_{\set{n+1} = P_1 \sqcup 
\dots \sqcup (P_s \cup \{n+1\}) \sqcup \dots \sqcup P_k} (p_1 - 1)! \dots 
p_s ! \dots (p_k-1)! \\
&= \sum_{s=1}^k (-1)^{k-1} \sum_{s=1}^k p_s \sum_{\set{n} = P_1 
\sqcup \dots \sqcup P_s \sqcup \dots \sqcup P_k} (p_1 - 1)! \dots (p_s - 
1)! \dots (p_k-1)! \\
&= n\sum_{s=1}^k (-1)^{k-1} \sum_{\set{n} = P_1 \sqcup \dots \sqcup 
P_k} (p_1 - 1)! \dots (p_k - 1)! \\
&= 0.
 \end{align*} 
\end{proof}
 
 \begin{proof}[of Theorem \ref{Th:ExtAct}]
Notice first that for any graph $F$ the statistical sum of {\em all} 
subgraphs $H \subset F$ equals 
 \begin{equation*}
Z(\{H \mid H \subset F\}) = \prod_{\edge{i}{j}\text{ is an edge of } F} 
(1+w_{ij}).
 \end{equation*}
Consider the poset $\mathcal P_n$ of all set partitions of 
$\set{n}$ ordered by refinement. In particular, $\name{min} = \{1\} \sqcup 
\dots \sqcup \{n\}$ is the smallest element of $\mathcal P_n$, and 
$\name{Max} = \set{n}$ is its largest element.

Recall (from Section \ref{Sec:Orient}) that for a set $P \subset \set{n}$ 
one denotes by  $\langle P\rangle$ the subgraph of $G$ with the vertex set 
$P$; edges of $\langle P\rangle$ are all the edges of $G$ with both 
endpoints in $P$. Statement \ref{St:Connect} implies now that 
 \begin{equation*}
Z(\{H \mid H \subset F\}) = \sum_{(P_1 \sqcup \dots \sqcup P_k) \in 
{\mathcal P}_n} C_{\langle P_1\rangle}(w) \dots C_{\langle P_k\rangle}(w). 
 \end{equation*}
By Lemma \ref{Lm:Moebius} the value $\mu(u,\name{Max}) = (-1)^{k-1} (k-1)!$ 
where $\mu(u,v)$ is the M\"obius function for the poset $\mathcal P_n$. 
Therefore theorem follows from the M\"obius inversion formula, see e.g.\ 
\cite{Stanley}.
 \end{proof}
 
\section*{Questions and final remarks}

Multivariate Tutte polynomial has been studied intensively since 1970s when 
it was found to be related to partition functions of some important models 
in mathematical physics (Ising model, Potts model, and more; for more 
information consult \cite{FK,Triang}, references therein and also the 
review paper \cite{Sokal}.) Of particular interest are the complex zeros of 
these polynomials because they are responsible for the phase transition in 
ferromagnetic and antiferromagnetic media. It might be very interesting to 
study the zeros of the external activity polynomial $G_G$. Many natural 
questions about them (including the half-plane property, see \cite{Sokal}) 
are still open. Notice however that by Corollary \ref{Cr:FreeTerm} the 
polynomial $C_G$ is related to the specialization of the Tutte polynomial 
at $q = 0$ while in Potts model $q$ is interpreted as a number of states 
(of the spin) --- we are leaving it to professional physicists to give a 
sensible interpretation to the polynomial $C_G$.

Another possible direction of study is suggested by the nature of formulas 
\eqref{Eq:MSub}, \eqref{Eq:MSubgrInv}, \eqref{Eq:GivenCycle}, 
\eqref{Eq:AltSum}, \eqref{Eq:SumD}, \eqref{Eq:TotCyclChr}, 
\eqref{Eq:TotCyclChi} and \eqref{Eq:ExtAct}: they all contain a  sign 
alternating summation. It is highly probable that these formulas present 
the Euler characteristics of suitable complexes; so there is a problem to 
find these complexes. They must be related to the categorification of the 
Tutte polynomial obtained recently by E.F.\,Jasso-Hernandez and Y.\,Rong in 
\cite{CategTutte}; see also \cite{Categ,Stosic1,Stosic2} where the 
categorification of the chromatic polynomial was carried out.

Section \ref{Sec:ChiZero} of this paper contains two different descriptions 
of graphs $G$ such that every connected component of $G$ has vanishing 
Euler characteristics. One is given by the Statement \ref{St:1Cycle} and 
its corollary (both are special cases of Theorem \ref{Th:2Core}), and the 
other is contained in Theorem \ref{Th:GenFun} and its corollary, based on 
the all-minors version of the matrix-tree theorem. The relation between 
these results resembles the relation between a determinant and its minors 
decomposition. It seems very interesting to find similar results for graphs 
with an arbitrary $2$-core.

\end{document}